\documentclass[11pt]{article}
\usepackage{amsmath}
\usepackage{amsfonts}
\usepackage{amssymb}

 \newtheorem{theorem}{Theorem}[section]

\def\be{\begin{equation}}
\def\ee{\end{equation}}
\def\bea{\begin{eqnarray}}
\def\eea{\end{eqnarray}}

\def\scri{{\cal{I}}}

\begin{document}
\title{Some examples of the behaviour of conformal geodesics}
\author{Paul Tod\\Mathematical Institute\\and\\St John's College,\\
Oxford}
\date{}
\maketitle
\begin{abstract}
With the aid of concrete examples, we consider the question of whether, in the presence of conformal curvature, a conformal geodesic can become trapped in smaller and smaller sets, or phrased informally: are spirals possible? 
We do not arrive at a definitive answer, but we are able to find situations where this behaviour is ruled out, including a reduction of the conformal-geodesic equations to quadratures in a specific non-conformally flat metric.
\end{abstract}
Keywords: conformal geodesics, conformal geometry, spirals

\noindent MSC: 53A30

\noindent Journal of Geometry and Physics subject classification: real and complex differential geometry
\section{Introduction}

Conformal geodesics are a family of curves naturally defined in a conformal manifold of any signature. 
They play a role with respect to the conformal structure similar to that played by the geodesics of a metric, and the definition has been around at least since the 1950s, \cite{KY}, but much less is known about them than 
about metric geodesics. For example there are unanswered questions related to completeness of the form: can a conformal geodesic enter and not leave smaller 
and smaller sets? A specific version of this question due to Helmut Friedrich \cite{HF} is: if a conformal geodesic $\Gamma$ enters every neighbourhood of a point $p$, does $\Gamma$ 
necessarily pass through $p$ (with a finite limiting velocity and acceleration)? There are similar and related open questions in Schmidt's construction of the $c$-boundary \cite{Sch}. 
A metric geodesic cannot be trapped in smaller and smaller sets 
and a proof can be based on the existence of geodesically-convex neighbourhoods: if a geodesic enters a geodesically-convex ball then it must leave it too. 
Without a similar notion of convexity for conformal 
geodesics, it is not clear what the answer to the corresponding question should be. 
In conformally-flat metrics, such as Euclidean or Minkowski spaces, there are conformal geodesics which are closed plane circles of arbitrarily 
small radius. It is usually thought that these will fail to close up if some conformal curvature is added (and we shall see examples confirming this below), but whether the perturbed 
conformal geodesics can then spiral in in some fashion, either to an arbitrarily small circle or to a fixed point, is not known.

In this article we approach the question of spirals with the aid of some examples, including a reduction to quadratures of the conformal geodesic equation in the Riemannian or 
Lorentzian metrics on \emph{Nil}, one of Thurston's geometries (see e.g. \cite{PS}; I believe this is the first published instance of a reduction to quadratures of the conformal geodesic equation in 
a non-conformally flat metric but, possibly just because this case is integrable, there are no spirals). We do not solve the problem of the existence of spirals in general 
but we are able to establish nonexistence in certain cases, and we find various 
other properties of conformal geodesics.

\medskip

The object of study then is a conformal manifold $(M,[g])$, where $[g]$ is a conformal equivalence class of nondegenerate, say smooth, metrics which can have 
any signature, though we shall usually restrict to Riemannian or Lorentzian (where this is taken to mean one negative eigenvalue for the metric).

Make a choice of metric $g_{ab}$ in the conformal class, with corresponding metric covariant
derivative $\nabla_a$. Then a conformal geodesic is a curve $\Gamma$ with
tangent vector $v^a$ and a one-form $b_a$ given along $\Gamma$ and
satisfying the system:
\bea v^c\nabla_cv^a&=&-2(v^cb_c)v^a+g^{ac}b_c(g_{ef}v^ev^f)  \label{cg1}\\
v^c\nabla_cb_a&=&(v^cb_c)b_a-\frac12g_{ac}v^c(g^{ef}b_eb_f)+L_{ac}v^c
\label{cg2}
 \eea
where, in dimension $n$, $L_{ab}=\frac{1}{(n-2)}(R_{ab}-\frac{1}{2(n-1)}Rg_{ab})$, $R_{ab}$ is the Ricci tensor and $R$ is the scalar curvature. Other equivalent definitions of conformal geodesic are possible (see \cite{BE}). In indefinite signatures, the definition 
still makes sense for null $v^a$ when (\ref{cg1}) reduces to the null geodesic equation: null geodesics are conformal geodesics, and the only 
null conformal geodesics are null geodesics. Henceforth we shall assume that the conformal geodesics considered are not null.

The conformal geodesic equations are conformally-invariant in that, under conformal
rescaling
\be\label{cs1}g_{ab}\rightarrow \tilde{g}_{ab}=\Omega^2g_{ab},\ee
with the usual corresponding change to $\nabla_a$, they take the same form
with
\bea \tilde{v}^a&=&v^a\label{cs2}\\
\tilde{b}_a&=&b_a-\Upsilon_a\label{cs3} \eea
where, as usual, $\Upsilon_a=\Omega^{-1}\nabla_a\Omega$.

A conformal geodesic admits a preferred parameter $\tau$ defined up
to choice of origin by
\be\label{cg3}v^c\nabla_c\tau=1.\ee
There is a reparametrisation freedom of M\"obius transformations in $\tau$:
\be\label{cg4} \tau\rightarrow\hat\tau=\frac{a\tau+b}{c\tau+d},\ee
so that $\tau$ may be called a projective parameter, which must be accompanied by a transformation of $v^a$ and $b_a$:
\bea\hat{v}^a&=&(c\tau+d)^{2}v^a\label{cg5}\\
\hat{b}_a&=&b_a+fg_{ac}v^c\label{cg6}\eea
where
\[f=-2(g_{ef}v^ev^f)^{-1}c(c\tau+d)^{-1}.\]
This transformation draws attention to a familiar problem in the study of conformal geodesics: it is possible for $\hat{v}^a$
to vanish and $\hat{b}_a$ to be singular at a regular point both of the manifold and of the curve, by choice of
projective parameter. One way to avoid this problem is to introduce a \emph{third-order form} of the equations.  How this is done depends slightly on the signature. In the Riemannian 
case or the Lorentzian case with a space-like velocity, introduce the unit tangent $u^a=\chi^{-2}v^a$, where
$\chi^4=g_{ef}v^ev^f>0$, then (\ref{cg1}) can be solved for $b_c$ in
terms of the acceleration $a^c=u^a\nabla_au^c$ as
\be\label{3a}b^c=a^c+(u^ab_a)u^c,\ee
where $u^ab_a=-2\chi^{-1}u^a\nabla_a\chi$. Use this in (\ref{cg2})
to find
\be\label{3b}u^c\nabla_ca^b=u^b(-g_{ef}a^ea^f-L_{ef}u^eu^f)+L^b_cu^c,\ee
which is the third-order form of the conformal geodesic equation, and
\[u^c\nabla_c(u^ab_a)=\frac12(g_{ef}a^ea^f+(u^ab_a)^2)+L_{ef}u^eu^f \]
which can be written as
\be\label{3c}\ddot{\chi}=-\frac14(g_{ef}a^ea^f+2L_{ef}u^eu^f)\chi,\ee
where the overdot denotes differentiation with respect to proper
length. We shall always call the parameter $t$, even when it is proper length. From (\ref{cg3}) and the definition of $u^a$ we have
\[\partial_\tau=v^a\partial_a=\chi^2u^a\partial_a=\chi^2\partial_t,\]
so that $d\tau=dt/\chi^2$.

Having solved (\ref{3b}), to recover a representative $(v^a,b_a)$ we
solve (\ref{3c}) for a choice of $\chi$ when $v^a=\chi^2u^a$, and
$b_a$ is obtained from (\ref{3a}). The different choices of $\chi$ correspond to the different choices of the projective  parameter $\tau$. 
For example, if we choose two solutions $\chi_1,\chi_2$ of (\ref{3c}) with unit Wronskian:
\[\chi_1\dot{\chi}_2-\chi_2\dot{\chi}_1=1,\mbox{   so that   }\frac{d}{dt}\left(\frac{\chi_2}{\chi_1}\right)=\frac{1}{\chi_1^2}\mbox{  and  }\chi_1^{-2}dt=d\left(\frac{\chi_2}{\chi_1}\right),\]
so that an allowed choice for $\tau$ is $\chi_2/\chi_1$. Thus zeroes and poles of the projective parameter are tied to zeroes and poles of the solutions 
of (\ref{3c}).

For the time-like Lorentzian case, introduce unit tangent $u^a=\chi^{-2}v^a$, where now
$\chi^4=-g_{ef}v^ev^f$, then (\ref{cg1}) can be solved for $b_c$ in
terms of the acceleration $a^c=u^a\nabla_au^c$ as
\be\label{3d}b^c=-a^c-(u^ab_a)u^c,\ee
where $u^ab_a=-2\chi^{-1}u^a\nabla_a\chi$.
Use this in (\ref{cg2})
to find
\be\label{3e}u^c\nabla_ca^b=u^b(g_{ef}a^ea^f-L_{ef}u^eu^f)-L^b_cu^c,\ee
and
\be\label{3f}\ddot{\chi}=\frac14(g_{ef}a^ea^f-2L_{ef}u^eu^f)\chi.\ee
Now for existence we have:
\begin{theorem} A conformal geodesic, that is a solution of (\ref{3b}) or (\ref{3e}), will always exist for some interval $(-t_1,t_1)$ in $t$, given data at time $t=0$ consisting of 
an initial position, an initial unit velocity $u^a$ and an initial acceleration $a^a$ orthogonal to $u^a$. Thus a conformal 
geodesic can always be continued in the Riemannian case as long as $a^a$ remains finite, and in the Lorentzian case as 
long as both $u^a$ and $a^a$ remain finite.
\end{theorem}
{\bf Proof} 

This follows at once from Picard's Theorem.    \begin{flushright} $\Box$ \end{flushright}

It follows as a corollary that spiralling can only occur in the Riemannian case when the acceleration diverges or in the Lorentzian case when either the velocity or the acceleration diverges (though presumably 
one cannot happen without the other). For conformal geodesics in flat space or in an Einstein manifold, the magnitude of the acceleration is constant. Thus spiralling cannot happen in a Riemannian Einstein 
space, and it seems likely that this is also true in Lorentzian Einstein spaces. In the presence of more general curvature, the acceleration need not be constant since we have, from (\ref{3b}) or (\ref{3e}), that
\[\frac{d}{dt}(g_{ef}a^ea^f)=\pm 2a^eu^fL_{ef}\]
with the sign depending on the signature. If we set $a^e=qw^e$ where $w^e$ is a unit vector and $q\geq 0$, then
\[\dot{q}=\pm w^eu^fL_{ef},\]
where the plus sign is for Riemannian or Lorentzian space-like and the minus sign is for Lorentzian time-like. Now for Riemannian, in a region in which the curvature 
is bounded, the magnitude $q$ of the acceleration has bounded time-derivative and so the acceleration can only diverge after infinite $t$, that is after infinite 
proper distance, and so spiralling could only occur after infinite proper distance. Furthermore, if the magnitude of 
the acceleration does diverge then by a simple comparison argument applied to (\ref{3c}), any admissible $\chi$ will have infinitely many zeroes, 
and so any choice of projective parameter $\tau$ will pass infinitely often through zero and infinity.

The Lorentzian cases are different, for example both the velocity and the acceleration can 
become unbounded in finite time, as we shall see from examples below.

\medskip

In the rest of this article, we briefly review in Section 2 the behaviour of conformal geodesics in three-dimensional Einstein spaces. Much of this carries over to conformally-flat 
spaces of any dimension and is quite untypical of the general case. 

In Section 3, we consider some very simple axisymmetric metrics in which it is possible to solve enough of the 
conformal geodesic equations to show that there are some which close up into circles but plenty more which do not: the widespread intuition that conformal curvature typically prevents 
circles from closing is vindicated. What makes the calculation tractable is that planes of constant $z$ in these metrics are totally-geodesic for conformal geodesics (call this property 
\emph{totally-conformally geodesic}). It is also possible to show in these examples that the conformal geodesics in these 2-dimensional totally-conformally geodesic surfaces cannot spiral 
into the origin or into arbitrarily small limit cycles round the origin - certain kinds of spiralling are ruled out. 

In Section 4, we reduce the conformal geodesic equations in a particular metric to quadratures. The metric is a Riemannian 
or Lorentzian version of \emph{Nil} (which arises in Thurston's classification of geometries, see e.g. \cite{PS}). The reduction is possible because of the existence of two constants of the motion. One comes 
from a Killing vector with a particular algebraic relation to the Ricci tensor; the other arises because the Ricci tensor is a Killing tensor. One of the constants serves to bound the norm-squared of 
the acceleration so that, in the Riemannian case, the acceleration cannot blow up and any conformal geodesic can be continued forever - there are no spirals in this example. In the Lorentzian case the indefiniteness of 
the signature means that a bound on the norm-squared of the acceleration does not give a bound on the components of acceleration and one can find explicit examples with the acceleration 
blowing up exponentially in time (which can already happen in Minkowski space) or in finite time (which cannot). Again, there are no spirals. 

In Appendix A we show 
that the conformal geodesic equations in the Berger sphere can also be reduced to quadratures. This is not too surprising given the result for \emph{Nil} because the metrics have a number of 
similarities, including the Killing tensor property of the Ricci tensor. 

In Appendix B we write out the conformal geodesic equation in the Schwarzschild metric, a Lorentzian Einstein metric. The 
novelty here is that there are two constants of the motion arising from the Killing spinor, though none from the Killing vectors so that we do not have enough constants for integrability.

\section{The cases of $S^3,E^3,H^3,M^3$}\label{sCF}
By this we mean the three, three-dimensional Riemannian Einstein metrics of positive, zero and negative curvature respectively, and three-dimensional Minkowski space. 
These cases are well understood but it is useful to review them for what follows. All are conformally-flat, so that the conformal 
geodesics in $S^3$ and $H^3$ can be obtained from those in $E^3$, and all are Einstein, 
so that the term in $L_{ab}$ drops 
out of (\ref{3b}) or (\ref{3e}). We're left with
\be\label{3g}u^c\nabla_ca^b=-\epsilon u^b(g_{ef}a^ea^f),\ee
where $\epsilon=-1$ for Lorentzian time-like and $\epsilon=1$ for Riemannian or
 Lorentzian space-like (so that $\epsilon=g_{ab}u^au^b$). We see at once that $g_{ef}a^ea^f$ 
is constant. Since this can be zero, we see that all metric geodesics are actually conformal geodesics. 
This will not be true with more general curvature: 

\medskip

\noindent{\bf Remark}: By inspection of (\ref{3b}) or (\ref{3e}), we see that only a metric geodesic whose tangent 
vector is an eigenvector of the Ricci tensor will also be a conformal geodesic.

\medskip

To find more constants of the motion, suppose $K^a$ is a Killing vector, so that
\[M_{ab}:=\nabla_aK_b=\nabla_{[a}K_{b]},\;\;\nabla_aM_{bc}=\frac16R(g_{ab}K_c-g_{ac}K_b),\]
where $R$ is the Ricci scalar, and we've used the expression
\[R_{abcd}=\frac16R(g_{ac}g_{bd}-g_{ad}g_{bc}),\]
for the Riemann tensor, which is valid in these metrics. Now consider the quantity
\be\label{3h}
Q:=M_{ab}u^aa^b-\frac16 RK_au^a.\ee
It is easy to see that this is constant, and this observation yields six constants of the motion, plus $g_{ef}a^ea^f$ and 
$g_{ef}u^eu^f$ for $S^3$ and $H^3$. For $E^3$ and $M^3$, where $R=0$, only the 3 rotations give non-trivial constants 
like this (the translations give zero), and these correspond to the constant bivector $u_{[a}a_{b]}$ which defines 
the plane of the orbit. 

To find the explicit solution for the conformal geodesics in $E^3$, when $\epsilon=1$, introduce $A^2=g_{ef}a^ea^f$ then 
(\ref{3g}) becomes
\[\ddot{u}^a=-A^2u^a.\]
If $A=0$ then the solution is a metric geodesic. If $A\neq 0$ then this is solved by
\[u^a=u^a_0\cos(At)+A^{-1}a^a_0\sin(At),\]
where $u^a_0=u^a(0)$ and $a^a_0=a^a(0)$; the position $x^a(t)$ follows by
integration:
\be\label{x1}
x^a=x^a_0+A^{-1}u^a_0\sin(At)+A^{-2}a^a_0(1-\cos(At)).\ee
The conformal geodesics are planar circles and straight lines.

To find the explicit solution in $M^3$, in the time-like case (so $\epsilon=-1$) $a^a$ is orthogonal to
$u^a$ and so is space-like. Now $g_{ef}a^ea^f=A^2$ and
(\ref{3g}) is solved by
\[u^a=u^a_0\cosh(At)+A^{-1}a^a_0\sinh(At),\]
whence again $x^a(t)$ follows by integration. This is the constant-acceleration world-line familiar from courses on special relativity: the scalars 
$g_{ab}u^au^b$ and $g_{ab}a^aa^b$ are constant but the components of both velocity and acceleration are unbounded in a constant basis.

The Lorentzian space-like case can have $a^a$ time-like,
space-like  or null. The assumption of time-like leads to another hyperbola, and space-like gives another circle, so the only new case is null, when $a^a=a^a_0$ is a
constant null vector and
\[x^a=x^a_0+tu^a_0+\frac12t^2a^a_0.\]
This is a parabola lying in the null hyperplane $g_{ab}a^a_0(x^b-x^b_0)=0$.

Solutions in $S^3$ or $H^3$ can be obtained from those in $E^3$ by conformal rescaling. They turn out to be intersections of $S^3$ embedded in the standard way in $E^4$ or $H^3$ embedded 
in the standard way in $M^4$ by linear subspaces in the ambient space, 
so they are circles or hyperbolae.
\section{An axisymmetric example}
We can get enough information to throw some light on our questions from a simple axisymmetric example. First note:

\begin{theorem}\label{thm1} Any trajectory of a (non-null) hypersurface-orthogonal Killing vector is a conformal geodesic.

\end{theorem}
{\bf Proof} 

The condition of hypersurface-orthogonality means
\[K_{[a}\nabla_bK_{c]}=0,\]
so that, since $K^a$ is a Killing vector,
\[\nabla_aK_b=\nabla_{[a}K_{b]}=2V_{[a}K_{b]},\]
for some vector $V_a$ which can be assumed to be orthogonal to $K^a$.

Suppose $\gamma$ is a trajectory of the Killing vector and set $K^aK_a=\epsilon U^2$ where $\epsilon=\pm 1$, to allow for different signatures and $U\neq 0$ by assumption.
 The unit tangent vector to $\gamma$ is $u^a=U^{-1}K^a$ and the acceleration is
\[a^a=u^b\nabla_bu^a=-\epsilon V^a,\]
using what we know. Recall the Killing vector identity
\[\nabla_a\nabla_bK_c=R_{bcad}K^d,\]
so that
\[R_{cd}K^d=K_c(\nabla_aV^a+V_aV^a)-K^b\nabla_bV_c.\]
and it is easy to check that  (\ref{3b}) is satisfied.   \begin{flushright} $\Box$ \end{flushright}

\noindent{\bf Remark:} There's no reason to expect this result to hold for Killing vectors which are not hypersurface-orthogonal, and examples where it fails will be seen below, 
in Section 4.1.

\medskip

Now we'll do the example. Suppose the metric is
\be\label{met1}
ds^2=dr^2+F^2(r)d\phi^2+dz^2,\ee
with a regular axis at $r=0$, and consider conformal geodesics lying in $z=0$. Introduce $\chi$ by
\[\dot{r}=\cos\chi,\;\;\;F\dot{\phi}=\sin\chi\]
to take care of the normalisation of the velocity (since $\dot{z}=0$), then the velocity is
\[\bf{u}=\cos\chi{\bf{e_1}}+\sin\chi{\bf{e_2}}\]
in terms of orthonormal vectors
\[{\bf{e_1}}=\partial_r,\;\;\;{\bf{e_2}}=F^{-1}\partial_\phi.\]
We claim the acceleration is then
\[{\bf a}=q(-\sin\chi{\bf{e_1}}+\cos\chi{\bf{e_2}}),\]
where
\[q=\dot{\chi}+\sin\chi F'/F,\]
and the prime means $d/dr$. Because of the dimension of the surface $z=0$, the Ricci term drops out of the conformal geodesic equation (\ref{3b}) and we're back to 
(\ref{3g}). The remaining equation is equivalent to the constancy of $g_{ab}a^aa^b$ or equivalently the constancy of $q$, which may be written
\[q=\dot{\chi}+\tan\chi\dot{ F}/F=\mbox{ constant},\]
where we've used $\dot{r}=\cos\chi$. Introduce
\[G(r)=\int_0^rF(s)ds,\]
then this integrates as
\[F\sin\chi-qG=C=\mbox{ constant},\]
from which we obtain
\[\dot{r}^2=1-\frac{(C+qG)^2}{F^2}=:f(r).\]
For a regular axis at $r=0$ we need $F/r\rightarrow 1$ as $r\rightarrow 0$, so also $2G/r^2\rightarrow 1$. Thus if $q$ and $C$ are both nonzero then $f(r)$ is negative for small $r$, 
and there is a lower limit for $r$ on the conformal geodesic, say $\Gamma$. Behaviour at large $r$ depends on what conditions we impose there. For asymptotic flatness 
we would want the same limits at large $r$ as at small $r$ and then, if $q$ and $C$ are both nonzero, $f(r)$ is negative for large $r$, so that 
$\Gamma$ is confined between finite and nonzero upper and lower limits in $r$. We also have the orbit equation 
as
 \[\phi=\int\frac{C+qG}{F(F^2-(C+qG)^2)^{1/2}}dr,\]
so that the conformal geodesics do not close in general, but rather fill out an annulus in the plane. This is a confirmation of the common belief noted in the Introduction.

The special case of $q=0$ is actually a metric geodesic, which can escape to infinity in $r$. The special case of $C=0$ is a conformal geodesic passing through the origin. If $q\neq 0$, it 
will have an upper limit in radius $r$ to which it will return after each pass through the origin - this is not a spiral but something more like a rosette. Indeed careful 
choices will give conformal geodesics which are closed rosettes, like the curve $r=a\sin 2\phi$.

From Theorem \ref{thm1}, it must be possible to choose $q$ and $C$ so that any fixed value of $r$, say $r=a$, is a conformal geodesic. This requires $f(a)=0=f'(a)$ 
which in turn requires
\[q=\pm F'(a)/F(a),\;\;C=\pm\left((F(a))^2-F'(a)G(a)\right)/F(a).\]
Now one might wonder whether there could be conformal geodesics which spiralled in asymptotically in time to these circles. The existence of such conformal geodesics for arbitrarily small $a$ 
would provide a negative answer to the question of Friedrich mentioned in the Introduction.

\medskip

\noindent{\bf Remark:} While this may be possible for some $a$, it cannot be arranged for arbitrarily small $a$. 

\medskip

To see 
this calculate
\[f''(a)=2\left(F(a)F''(a)-(F'(a))^2\right)/(F(a))^2,\]
at a circular orbit. At small enough $a$ this must be negative since $F\rightarrow 0,F'\rightarrow 1$ as $a\rightarrow 0$. Thus small 
circular conformal geodesics are stable, but to be the limit cycle of a spiral, a small circle 
would need to be unstable (or at worst marginally stable).

\medskip

Spirals can also be ruled out in a slightly more general case: replace (\ref{met1}) by
\be\label{met3}
ds^2=dr^2+F^2(r,z)d\phi^2+H^2(r,z)dz^2,
\ee
but where we'll suppose $F$ and $H$ are even in $z$ so that the plane $z=0$ is totally-geodesic for metric geodesics. We shall see below that conformal geodesics can also be confined to 
this plane. The velocity is
\[\alpha=\dot{r},\;\;\beta=F\dot{\phi},\;\;\mbox{  so that  }\alpha^2+\beta^2=1.\]
Then the acceleration is given by
\[a_1=\dot{\alpha}-\beta^2\frac{F'}{F},\;\;a_2=\dot{\beta}+\alpha\beta\frac{F'}{F},\]
retaining prime for $\partial_r$, and the conformal geodesic equations turn out to be
\[\dot{a}_1=\beta a_2\frac{F'}{F}-|a|^2\alpha+\alpha\beta^2Q,\]
\[\dot{a}_2=-\beta a_1\frac{F'}{F}-|a|^2\beta-\alpha^2\beta Q,\]
where
\[Q=-\frac{H''}{H}+\frac{F_{zz}}{FH^2}+\frac{F'H'}{FH},\]
which is a combination of curvature components, and in particular must be bounded as $r\rightarrow 0$. 

Since the velocity is a unit vector we must have
\[a_1\alpha+a_2\beta=0\mbox{   so   }a_1=q\beta, \;\;a_2=-q\alpha,\]
for some $q$, when in fact $q^2=|a|^2$. Now the conformal geodesic equations collapse to
\[\dot{q}=\alpha\beta Q=\dot{r}\beta Q,\]
so $dq/dr=\beta Q$ which is bounded for small $r$. Therefore one cannot have a conformal geodesic spiralling in to the origin since along such a curve the magnitude of the acceleration must be unbounded, 
while this equation would force $q$ to be bounded.

It is also possible to check in this more complicated example that while the axial Killing vector trajectories provide small circular conformal geodesics, 
the smallest circles are again stable and cannot be the future limit cycles of spirals.

\medskip

We've just seen that the plane $z=0$ is totally-conformally-geodesic in that conformal geodesics could be confined to it. It is worth noting the following:

\begin{theorem} The necessary and sufficient condition for a hypersurface $\Sigma$ to be totally-conformally-geodesic is the vanishing of the 
second-fundamental form of $\Sigma$.
\end{theorem}

{\bf Proof} 

Let the unit normal to $\Sigma$ be $N^a$, then the covariant derivative of $N_a$ takes the form
\[\nabla_aN_b=N_aA_b+K_{ab},\]
where
\[N^aA_a=0=K_{ab}N^b,\;\;K_{[ab]}=0,\]
so $K_{ab}$ is the second fundamental form of $\Sigma$. Start a conformal geodesic in $\Sigma$ so that initially 
\[q_1:=N_au^a=0,\;\;q_2:=N_aa^a=0.\]
Calculate
\[\dot{q}_1=q_2+K_{ab}u^au^b+q_1(u^aA_a),\]
and
\[\dot{q_2}=-q_1(g_{ab}a^aa^b)+N^aL_{ab}u^b+u^aa^b(K_{ab}+N_aA_b).\]
Now uniqueness of solution of this system will force $q_1=0=q_2$ for every initial choice provided $K_{ab}u^au^b=0=N^aL_{ab}u^b$ for every allowed $u^a$. This certainly requires $K_{ab}=0$. Suppose 
that holds and recall the Gauss-Codazzi equation:
\[R_{bc}h^c_aN^b=D_bK^b_a-D_aK,\]
where $h_{ab}=g_{ab}-N_aN_b$ is the intrinsic metric of $\Sigma$ and $D_a$ is the metric covariant derivative of $h$. Thus the vanishing of $K_{ab}$ 
implies the vanishing of $N^aL_{ab}h^b_c$, which gives the second condition. \begin{flushright} $\Box$ \end{flushright}

\section{Conformal geodesics in \emph{Nil}}
\emph{Nil} is one of Thurston's list of eight geometries \cite{PS}. It admits a transitive 4-dimensional isometry group isomorphic to the 4-dimensional nilpotent Lie group. 
It can be considered with different signatures and since, as we shall see, the conformal geodesic equation can be reduced to quadratures, it gives information on the behaviour 
of conformal geodesics in the non-conformally-flat setting.

We may take the metric to be
\be\label{n1}
g=dx^2+\epsilon dy^2+(dz-\lambda xdy)^2,\ee
with $\epsilon =\pm 1$ in order to consider either Riemannian or Lorentzian forms.
\subsection{Riemannian \emph{Nil}}
We'll start with Riemannian, so $\epsilon=1$, and choose an orthonormal basis of one-forms:
\be\label{n2}
\theta^1=dx,\;\;\theta^2=dy,\;\;\theta^3=dz-\lambda xdy,\ee
with dual basis
\be\label{n3}
e_1=\partial_x,\;\;e_2=\partial_y+\lambda x\partial_z,\;\;e_3=\partial_z.\ee
We calculate the connection one-forms as
\[\omega^1_{\;\;\;2}=\frac{\lambda}{2}\theta^3,\;\;\omega^1_{\;\;\;3}=\frac{\lambda}{2}\theta^2,\;\;\omega^2_{\;\;\;3}=-\frac{\lambda}{2}\theta^1,\;\;\]
and the nonzero Ricci tensor components as
\[R_{11}=R_{22}=-R_{33}=-\lambda^2/2.\]
We expand the unit velocity of a conformal geodesic in the chosen basis as 
\[u=\alpha e_1+\beta e_2+\gamma e_3,\]
so that
\be\label{nn1}\dot{x}=\alpha,\;\;\dot{y}=\beta,\;\;\dot{z}-\lambda x\dot{y}=\gamma,\ee
and $\alpha^2+\beta^2+\gamma^2=1$. 
The acceleration is then calculated to be 
\[a=a_1e_1+a_2e_2+a_3e_3\]
where
\bea
a_1&=&\dot{\alpha}+\lambda\beta\gamma\\\nonumber
a_2&=&\dot{\beta}-\lambda\alpha\gamma\\\label{ac1}
a_3&=&\dot{\gamma}\nonumber
\eea
and the conformal geodesic equations then are
\bea\label{ac2}
\dot{a}_1&=&-\frac{\lambda}{2}\gamma a_2-\frac{\lambda}{2}\beta a_3-\alpha|a|^2-\lambda^2\alpha\gamma^2\\\nonumber
\dot{a}_2&=&\frac{\lambda}{2}\gamma a_1+\frac{\lambda}{2}\alpha a_3-\beta|a|^2-\lambda^2\beta\gamma^2\\
\dot{a}_3&=&\frac{\lambda}{2}\beta a_1-\frac{\lambda}{2}\alpha a_2-\gamma|a|^2+\lambda^2\gamma(\alpha^2+\beta^2)\nonumber
\eea
where $|a|^2=g_{ab}a^aa^b$. From these, we claim that
\be\label{n6}
E:=|a|^2-\lambda^2\gamma^2=\mbox{ constant.}\ee
This is just a matter of checking. To understand 
the constant $E$ note that in general
\[\frac{d}{dt}(g_{ab}a^aa^b)=2L_{ab}u^aa^b=\frac{d}{dt}(L_{ab}u^au^b)-u^au^bu^c\nabla_aL_{bc},\]
and for \emph{Nil} the Ricci tensor satisfies the \emph{third Ledger Condition} \cite{PT}, equivalently the Killing tensor equation:
\[\nabla_{(a}R_{bc)}=0,\]
so that $|a|^2-L_{ab}u^au^b$ is constant since $u^au^bu^c\nabla_aL_{bc}=0$. $E$ is a combination of this constant and a multiple of $|u|^2$.

\medskip

\noindent{\bf Remark}:  The existence of the constant $E$ provides a bound on $|a|^2$, since necessarily $\gamma^2\leq 1$. In turn this shows that, in Riemannian \emph{Nil}, a conformal geodesic can be continued forever: 
there will not be any spirals in this case. 

\medskip

Next note that the quantity 
\be\label{n7}
J:=a_1\beta-a_2\alpha+\frac12\lambda\gamma\ee
is another constant of the motion. This is again a simple matter of checking, and this constant arises from the Killing vector $\partial_z$ analogously to the constants obtained in (\ref{3h}) from Killing vectors. 

It will follow from the existence of these two constants of the motion, that the equations are integrable. Note from the third of (\ref{ac2}) that
\[\ddot{\gamma}=\dot{a_3}= \frac{\lambda}{2}(J-\frac12\lambda\gamma)-\gamma(E+\lambda^2\gamma^2)+\lambda^2\gamma(1-\gamma^2),\]
which integrates to give
\[\dot{\gamma}^2=-\lambda^2\gamma^4+(\frac34\lambda^2-E)\gamma^2+\lambda J\gamma +C_1,\]
in terms of a new constant $C_1$, which we'll fix shortly. Suppose $\gamma$ known and parametrise $\alpha,\beta$ by
\[\alpha=(1-\gamma^2)^{1/2}\cos\chi,\;\;\;\beta=(1-\gamma^2)^{1/2}\sin\chi.\]
then the definition of $J$ leads to
\be\label{n8}
\dot{\chi}=(-J+\frac32\lambda\gamma-\lambda\gamma^3)(1-\gamma^2)^{-1}.\ee
Now we express the acceleration in terms of $\chi,\gamma$ and $\dot{\gamma}$ to discover that $C_1=E-J^2$. 

\medskip

We have 
reduced the problem to quadratures: first for $\gamma$ we have
\be\label{n9}
\dot{\gamma}^2=-\lambda^2\gamma^4+(\frac34\lambda^2-E)\gamma^2+\lambda J\gamma+E-J^2=:F(\gamma);\ee
then $\chi$ is obtained from (\ref{n8}) and the coordinates from (\ref{nn1}). All quantities will be analytic in time, since the solution 
$\gamma$ of (\ref{n9}) is. Note in (\ref{n9}) that $F(\pm 1)=-(J\mp\frac{\lambda}{2})^2$ so that, as desired, $\gamma$ is confined within a range which can extend 
to $1$ if $J=\lambda/2$ or to $-1$ if $J=-\lambda/2$, but not both; for general $J$ neither limit is attained and for no orbit are both attained. $\gamma$ will be periodic in time but $\chi$ 
will not in general, so that $\alpha,\beta$ also will not. In particular the orbits will not close in general.

One example which can be solved explicitly is $E=0=J$. Then 
\[\dot{\gamma}^2=\lambda^2\gamma^2(\frac34-\gamma^2),\]
which is solved by
\[\gamma=\frac{\sqrt{3}}{2}\mbox{sech}\omega,\mbox{   with   }\dot{\omega}=\frac{\lambda\sqrt{3}}{2}.\]
We may choose the origin of time so that $\omega=\frac{\lambda\sqrt{3}}{2}t$ and solve for $\chi$ to find
\[\tan(\chi-\chi_\infty)=\frac{3\sinh\omega}{1-2\sinh^2\omega},\]
where $\chi_\infty$ is a constant of integration. Now as $t\rightarrow\infty$ we have $\gamma\rightarrow 0$ and $\chi$ tending to a constant. Thus $\alpha$ and $\beta$ will tend to constants and 
the acceleration will tend to zero: this class of conformal geodesics 
tend to metric geodesics. It is easy to see that the metric geodesics which are also conformal geodesics fall into two classes: those with $\alpha=\beta=0$ and $\gamma=\pm 1$, 
and those with $\gamma=0$ and $\alpha,\beta$ constants. It is the second set which are involved here.

\medskip

\noindent{\bf Remark:} We can use the general solution to justify the claim made above, that the integral curves of a Killing vector which is not hypersurface orthgonal 
need not be conformal geodesics. 

\medskip

To see this, consider the Killing vector
\[K=\beta\partial_y+\gamma\partial_z\]
for constant $\beta,\gamma$ (which will shortly have the usual meanings). As a one-form this is
\[K=\beta dy+(\gamma-\lambda\beta x)(dz-\lambda xdy)\]
so that $K\wedge dK=\lambda(\beta^2-(\gamma-\lambda\beta x)^2)dx\wedge dy\wedge dz\neq 0$. Any trajectory of this Killing vector with $x=0$ has unit tangent $u=K=\beta e_2+\gamma e_3$, 
assuming that $\beta^2+\gamma^2=1$, and has acceleration
\[a_1=\lambda\beta\gamma,\;\;a_2=a_3=0.\]
Substituting into (\ref{ac2}) we see that the trajectory is a conformal geodesic iff $\beta\gamma=0$. So if we choose a Killing vector with $\beta\gamma\neq 0$, then its integral curves are not 
conformal geodesics.

\subsection{Lorentzian \emph{Nil}}
For the Lorentzian case, we take the metric to be 
\[g=dx^2- dy^2+(dz-\lambda xdy)^2.\]
We'll just consider space-like curves, so retaining (\ref{nn1}) $\alpha^2-\beta^2+\gamma^2=1$. Then (\ref{ac1}) becomes
\bea
a_1&=&\dot{\alpha}+\lambda\beta\gamma\\\nonumber
a_2&=&\dot{\beta}+\lambda\alpha\gamma\\\label{ac11}
a_3&=&\dot{\gamma}\nonumber
\eea
and (\ref{ac2}) becomes
\bea
\dot{a}_1&=&-\frac{\lambda}{2}\gamma a_2-\frac{\lambda}{2}\beta a_3-\alpha(a_1^2-a_2^2+a_3^2)+\lambda^2\alpha\gamma^2\\\nonumber
\dot{a}_2&=&-\frac{\lambda}{2}\gamma a_1-\frac{\lambda}{2}\alpha a_3-\beta(a_1^2-a_2^2+a_3^2)+\lambda^2\beta\gamma^2\\\label{ac3}
\dot{a}_3&=&\frac{\lambda}{2}\beta a_1-\frac{\lambda}{2}\alpha a_2-\gamma(a_1^2-a_2^2+a_3^2)-\lambda^2\gamma(1-\gamma^2)\nonumber
\eea
The constants of the motion become
\be\label{ac4}E:=a_1^2-a_2^2+a_3^2+\lambda^2\gamma^2=\mbox{constant},\ee
which does not now give a bound on the $a_i$, and indeed the acceleration will not now be bounded in general, and
\[J:=a_1\beta-a_2\alpha-\frac12\lambda\gamma=\mbox{ constant.}\]
The system reduces to quadratures again but now with
\be\label{L9}
\dot{\gamma}^2=F(\gamma):=\lambda^2\gamma^4-(\frac34\lambda^2+E)\gamma^2+\lambda J\gamma+E+J^2.\ee
Note that $F(\pm 1)=(J\pm\lambda/2)^2$, so $\gamma$ can go through $\pm 1$. 

The principal novelty in the Lorentzian case is that $\gamma$ can blow up in finite time: from (\ref{L9}) once $\gamma$ 
is large we can have $\dot{\gamma}\sim\pm\lambda\gamma^2$ so $\gamma\sim\pm(\lambda(t_0-t))^{-1}$ for some $t_0$. 

To understand the conformal geodesics, we need to analyse the polynomial $F$, which is a quartic and may have 0, 2 or 4 zeroes. If it 
has none then $\gamma$ necessarily ranges from $-\infty$ to $+\infty$, or vice versa; if it has 2 then $\gamma$ comes in from $\pm\infty$, meets a zero of $F$ and goes back to 
$\pm\infty$. If $F$ has 4 zeroes then there are trajectories with unbounded $\gamma$ and others where $\gamma$ is confined between finite upper and lower limits, and so is periodic.

For $\alpha$ and $\beta$ we distinguish two cases:
\begin{itemize}
 \item if $\gamma^2>1$, set $\alpha=(\gamma^2-1)^{1/2}\sinh\chi,\;\;\beta=(\gamma^2-1)^{1/2}\cosh\chi$;
\item if $\gamma^2<1$, set $\alpha=(1-\gamma^2)^{1/2}\cosh\chi,\;\;\beta=(1-\gamma^2)^{1/2}\sinh\chi$.
\end{itemize}
In either case, from the definition of $J$ we obtain
\[\dot{\chi}=(J+\frac32\lambda\gamma-\lambda\gamma^3)(\gamma^2-1)^{-1}.\]
This equation is singular at $\gamma^2=1$ but this is not really a problem. An unbounded $\gamma$ necessarily entails an unbounded $\beta$, but it is not at this point clear how $\alpha$ behaves.  
A periodic $\gamma$ can give a nonzero secular change in $\chi$ over a period; if it does not then $\alpha,\beta$ are also periodic; otherwise they grow with one of 
$\alpha\pm\beta\rightarrow\infty$ i.e. the velocity is asymptotically null.

We again consider the special case $E=0=J$. Then (\ref{L9}) becomes
\[\dot{\gamma}^2=\lambda^2\gamma^2(\gamma^2-\frac34),\]
which integrates to give $\gamma =\frac{\sqrt{3}}{2}\sec\omega$, where again $\omega= \frac{\lambda\sqrt{3}}{2}t$ with suitable choice of origin in $t$. Evidently $\gamma\rightarrow\pm\infty$ 
as $t\rightarrow\pm\pi/(\lambda\sqrt{3})$: a finite time blow-up for the velocity. The integration for $\alpha$ and $\beta$ leads to
\[\alpha+\beta=C_1\frac{(1-2\sin\omega)}{2(1+\sin\omega)},\;\;\alpha-\beta=C_2\frac{(1+2\sin\omega)}{2(1-\sin\omega)},\]
where $C_1$ and $C_2$ are fixed by data at $t=0$, with $C_1C_2=1$. The combinations $\alpha\pm\beta$ blow up at $\omega=\mp\pi/2$ respectively. From 
(\ref{nn1}) it follows that the coordinate combinations $x\pm y$ will diverge in these limits too. It follows from (\ref{ac4}) that the components of acceleration and the scalar $g_{ef}a^ea^f$ are 
also blowing up in finite time, which as we've seen cannot happen in the Riemannian \emph{Nil} case. 

In any case in which $\gamma$ becomes large its clear from (\ref{L9}) that $\gamma$ runs away 
to infinity in finite time so that the terms in $E$ and $J$ in (\ref{L9}) become negligible and the case that we've just integrated can be seen to be typical. From the discussion of the zeroes of $F$ 
we can say that, in any case in which $\gamma$ becomes large, the coordinates $x$ and $y$ diverge so there are no spirals, while if $\gamma$ is periodic then either $\alpha,\beta$ are also periodic or 
one of $\alpha\pm\beta$ diverges, whence so does one of $x\pm y$ and there are no spirals in this case. 

It remains to consider solutions with constant $\gamma$. If $\gamma=\gamma_0=\mbox{const.}$ then the constants $E$ and $J$ and the normalisation of the velocity give
\begin{eqnarray*}
\alpha^2-\beta^2&=&1-\gamma_0^2\\
a_1\alpha-a_2\beta &=&0\\
a_1\beta-a_2\alpha&=&J+\frac12\lambda\gamma_0\\
a_1^2-a_2^2&=&E-\lambda^2\gamma^2
\end{eqnarray*}
whence $a_1=k\beta,a_2=k\alpha$ for constant $k=(J+\frac12\lambda\gamma_0)(\gamma_0^2-1)^{-1}$. Now from (\ref{ac11})
\[k\beta=a_1=\dot{\alpha}+\lambda\beta\gamma_0,\]
and similarly for $\beta$, so that
\[\dot{\alpha}=\Omega\beta,\;\;\dot{\beta}=\Omega\alpha,\]
with $\Omega=k-\lambda\gamma_0$. The solutions are exponentials (or constants if $\Omega=0$). Now the coordinates $x,y$ will eventually be monotonic in time, diverging 
to $\pm\infty$ and there are no spirals in this case.

An explicitly solvable special case with constant $\gamma$ is when $\alpha=\beta$; then $\gamma^2=1$ so we may suppose $\gamma=1$. The acceleration is
\[a_1=a_2=\dot{\alpha}+\lambda\alpha,\;\;a_3=0\]
so that the acceleration is a null vector. The conformal geodesic equations collapse to
\[\ddot{\alpha}+\lambda\dot{\alpha}=-\frac12\lambda(\dot{\alpha}+\lambda\alpha)+\lambda^2\alpha,\]
or
\[\ddot{\alpha}+\frac32\lambda\dot{\alpha}-\frac12\lambda^2\alpha=0.\]
The solution is a sum of exponentials so that the velocity and acceleration are both unbounded (though bounded at any finite time). The orbit is not a spiral.

\medskip

In conclusion, in the case of the Riemannian or Lorentzian metrics on \emph{Nil}, the conformal-geodesic equations can be reduced to quadratures and there are no spiralling 
conformal geodesics.

\section*{Acknowledgements}
I gratefully acknowledge useful discussions with Helmut Friedrich, Bernd Schmidt, Michael Eastwood and Robert Bryant and wish to thank the Albert 
Einstein Institute in Golm and the Mathematics Department of the Australian National University, where this work was begun.

\section*{Appendix A: Conformal geodesics in the Berger Sphere}
\setcounter{equation}{0}
\renewcommand{\theequation}{A\arabic{equation}}
We relegate to this Appendix the discussion of conformal geodesics in the metric known as the Berger sphere. This case also turns out to reduce to quadratures, 
but does not lead to new phenomena. The metric is 
\[g=(\sigma^1)^2+(\sigma^2)^2+\lambda^2(\sigma^3)^2,\]
where $\sigma^i$ are the usual left-invariant one-forms on $SU(2)$, so that
\[d\sigma^1=\sigma^2\wedge\sigma^3.\]
In Euler angles these one-forms can be taken to be
\[\sigma^1+i\sigma^2=e^{i\psi}(d\theta-i\sin\theta d\phi),\;\;\sigma^3=d\psi+\cos\theta d\phi,\]
when the dual basis of vector-fields is 
\[e^{i\psi}(\partial_\theta-\frac{i}{\sin\theta}\partial_\phi+i\cot\theta\partial_\psi),\;\;\partial_\psi.\]
Clearly the metric becomes the round metric on ${\bf S}^3$ when $\lambda=1$. For the curvature calculation, we take the orthonormal basis of forms to be
\[\theta^1=\sigma^1,\;\;\theta^2=\sigma^2,\;\;\theta^3=\lambda\sigma^3,\]
then the connection one-forms are
\[\omega^1_{\;\;2}=C\theta^3,\;\;\omega^2_{\;\;3}=A\theta^1,\;\;\omega^3_{\;\;1}=A\theta^2,\]
where $A=\lambda/2$ and $C=(2-\lambda^2)/2\lambda$. The nonzero Ricci components are
\[R_{11}=R_{22}=1-\lambda^2/2,\;\;R_{33}=\lambda^2/2,\]
so the curvature scalar is $R=2-\lambda^2/2$.

We expand the velocity $u$ and acceleration $a$ in terms of the basis $(e_1,e_2,e_3)$ of vector fields dual to the basis $(\theta^i)$ as
\[u=\alpha e_1+\beta e_2+\gamma e_3,\]
\[a=a_1e_1+a_2e_2+a_3e_3,\]
then the acceleration has components
\begin{eqnarray*}
a_1&=&\dot{\alpha}+(C-A)\beta\gamma\\
a_2&=&\dot{\beta}-(C-A)\gamma\alpha\\
a_3&=&\dot{\gamma},
\end{eqnarray*}
and the conformal geodesic equations become
\begin{eqnarray*}
\dot{a_1}&=&-C\gamma a_2+A\beta a_3-\alpha(a_1^2+a_2^2+a_3^2)+R_{11}\alpha(1-\alpha^2-\beta^2)-R_{33}\alpha\gamma^2\\ 
\dot{a_2}&=&-A\alpha a_3+C\gamma a_1-\beta(a_1^2+a_2^2+a_3^2)+R_{11}\beta(1-\alpha^2-\beta^2)-R_{33}\beta\gamma^2\\ 
\dot{a_3}&=&-A\beta a_1+A\alpha a_2-\gamma(a_1^2+a_2^2+a_3^2)-R_{11}\gamma(\alpha^2+\beta^2)+R_{33}\gamma(1-\gamma^2).
\end{eqnarray*}
As before we seek constants of the motion, and find
\[J:=a_1\beta-a_2\alpha-A\gamma=\mbox{  constant},\]
which is analogous to (\ref{n7}) and obtained from the Killing vector $\partial_\psi$, and
\[E:=|a|^2+\gamma^2(R_{11}-R_{33})=a_1^2+a_2^2+a_3^2+(1-\lambda^2)\gamma^2=\mbox{  constant},\]
which is analogous to (\ref{n6}). $E$ serves to bound the magnitude of the acceleration, and exists because the Ricci tensor is again a Killing tensor. 
We derive an equation for $\gamma$ alone:
\[\ddot{\gamma}=(1-\lambda^2)(2\gamma^3-\gamma)-(E+A^2)\gamma-AJ\]
We shall integrate this shortly, but before that set
\[\alpha=R\cos\chi,\;\;\beta=R\sin\chi,\mbox{  where  }R=(1-\gamma^2)^{1/2},\]
for then the definition of $J$ leads to
\[\dot{\chi}=-\frac{J+A\gamma-(C-A)\gamma(1-\gamma^2)}{1-\gamma^2}.\]
and then
\[(\dot{\gamma})^2=F(\gamma):=(1-\lambda^2)(\gamma^4-\gamma^2)-(E+A^2)\gamma^2-2AJ\gamma +E-J^2.\]
Note that $F(\pm 1)=-(J\pm A)^2$, so that the evolution forces $\gamma^2\leq 1$. The general discussion is now just like the case of 
Riemannian \emph{Nil}: since the acceleration is bounded, the orbits can be extended indefinitely but will not close in general. 

The special solutions with $E=0=J$ can be found explicitly:
\[(\dot{\gamma})^2=(1-\lambda^2)(\gamma^4-\gamma^2)-A^2\gamma^2=\gamma^2((1-\lambda^2)\gamma^2-(1-\frac34\lambda^2)).\]
These do not exist for  $4/3>\lambda^2$ but for  $\lambda^2>4/3$ we find

\[\gamma=k\,\mathrm{sech}\mu t,\;\;\mbox{  with }k^2=\frac{3\lambda^2/4-1}{\lambda^2-1},\;\;\mu=k(\lambda^2-1)^{1/2}.\]
Now $\gamma,\dot{\gamma}\rightarrow 0$ as $t\rightarrow \pm\infty$, so that in these limits the acceleration goes to zero and the conformal geodesic asymptotes to a metric geodesic. 
The metric geodesics which are conformal geodesics form two classes: those with $\gamma^2=1$ so that $\alpha=0=\beta$; and those with $\gamma=0$ and $\alpha$, $\beta$ constants. It is 
the second family which are serving as asymptotes here. Brief calculation shows this class are closed (as are the other class) but this is not true for all geodesics.

As with \emph{Nil} there is no spiralling, which may be a reflection of the integrability.

\section*{Appendix B: Conformal geodesics in the Schwarzschild solution}
\setcounter{equation}{0}
\renewcommand{\theequation}{B\arabic{equation}}
The Schwarzschild solution is a familiar Lorentzian Einstein metric, so that the conformal geodesic equations imply that the 
magnitude of the acceleration is constant, but there is no reason to expect the components of the acceleration to be bounded. 
This case is worth studying as there are two constants of the motion arising from the Killing spinor, although 
this is not sufficient to reduce the system to quadratures. 

The metric in the region outside the horizon can be taken to be
\be\label{s1}
-V^2dt^2+\frac{1}{V^2}dr^2+r^2(d\theta^2+\sin^2\theta d\phi^2),\ee
with $V^2=1-2m/r$, (so we are restricting to $r>2m$). We take the orthonormal basis of one-forms as
\[\theta^0=Vdt,\;\;\theta^1=V^{-1}dr,\;\;\theta^2=rd\theta,\;\;\theta^3=r\sin\theta d\phi.\]
Then expand in this basis:
\[u=\alpha\theta^0+\beta\theta^1+\gamma\theta^2+\delta\theta^3\]
\[a=a_0\theta^0+a_1\theta^1+a_2\theta^2+a_3\theta^3\]
and calculate
\begin{eqnarray*}
\dot{\alpha}&=&a_0-\alpha\beta V'\\
\dot{\beta}&=&a_1-\alpha^2 V'+(\gamma^2+\delta^2)V/r\\
\dot{\gamma}&=&a_2-\beta\gamma V/r+\delta^2\cot\theta/r\\
\dot{\delta}&=&a_3-\beta\delta V/r-\gamma\delta\cot\theta/r\\
\dot{a}_0&=&-\epsilon A\alpha-\alpha a_1V'\\
\dot{a}_1&=&-\epsilon A\beta-\alpha a_0V'+(\gamma a_2+\delta a_3)V/r\\
\dot{a}_2&=&-\epsilon A\gamma-\gamma a_1V/r+\delta a_3\cot\theta /r\\
\dot{a}_3&=&-\epsilon A\delta-\delta a_1V/r-\delta a_2\cot\theta /r
\end{eqnarray*}
where
\[A=-a_0^2+a_1^2+a_2^2+a_3^2.\]
This is the norm of the acceleration which is constant but may have either sign, and $\epsilon$ is -1 for time-like $u$ and +1 for space-like,
so that
\[\epsilon=-\alpha^2+\beta^2+\gamma^2+\delta^2,\]
and of course
\[-\alpha a_0+\beta a_1+\gamma a_2+\delta a_3=0.\] 
For time-like $u$, the acceleration is space-like and  $A$ must be positive. The following can be checked to be constants of the motion:
\[r(\beta a_0-\alpha a_1)-\epsilon\alpha V,\;\;r(\delta a_2-\gamma a_3).\]
These are the real and imaginary parts of the quantity
\be\label{s2}Q=u^aa^bF_{ab}-\epsilon K_bu^b,\ee
where $F_{ab}=\omega_{AB}\epsilon_{A'B'}$ is obtained from the Killing spinor $\omega_{AB}$ which satisfies
\[\nabla_{AA'}\omega_{BC}=2K_{A'(B}\epsilon_{C)A}\]
with $K_a=K_{AA'}$ the time-like Killing vector, unit at infinity. There do not seem to be constants of the motion related directly to the (four linearly-independent) Killing vectors 
of the solution so it is interesting that two arise from the Killing spinor. From the definition, it is simple to check that $Q$ in (\ref{s2}) is a constant of the motion, and will be as well for the Kerr solution, which 
likewise possesses a Killing spinor.

We may explore across the horizon by introducing $\lambda=V(\alpha+\beta)/2,\nu=(\alpha-\beta)/V$ and correspondingly $a_L=V(a_0+a_1)/2,a_N=(a_0-a_1)/V$, which follow by introducing 
the two null vectors $L=(e_0+e_1)/V,N=V(e_0-e_1)/2$ 
into the basis in place of $e_0,e_1$. The corresponding null coordinates are obtained by solving
\be\label{n4}
du=dt-\frac{dr}{V^2},\;\;\;dv=dt+\frac{dr}{V^2},
\ee
and the system regular at the horizon is 
\begin{eqnarray*}
\dot{\lambda}&=&a_L+\frac{V^2}{2r}(2\lambda\nu+\epsilon)-VV'\lambda\nu \\
\dot{\nu}&=&a_N-\frac1r(2\lambda\nu+\epsilon)+VV'\nu^2 \\
\dot{\gamma}&=&a_2-\frac{\gamma}{r}(\lambda-\frac12\nu V^2)+\delta^2\cot\theta/r\\
\dot{\delta}&=&a_3-\frac{\delta}{r} (\lambda-\frac12\nu V^2)-\gamma\delta\cot\theta/r\\
\dot{a}_L&=&- \epsilon A\lambda-VV'\nu a_L+\frac{V^2}{2r}(\nu a_L+\lambda a_N)\\
\dot{a}_N&=&-\epsilon A\nu+VV'\nu a_N-\frac1r(\nu a_L+\lambda a_N)\\
\dot{a}_2&=&-\epsilon A\gamma-\frac{\gamma}{r}(a_L-\frac12a_NV^2)+\delta a_3\cot\theta /r\\
\dot{a}_3&=&-\epsilon A\delta-\frac{\delta}{r}(a_L-\frac12a_NV^2)-\delta a_2\cot\theta /r\\
\end{eqnarray*}
 The first constant of the motion can be written
\[r(\lambda a_N-\nu a_L)-\epsilon (\lambda+\frac12V^2\nu).\]

Without integrability we are restricted to integrating some special cases. The case of radial time-like conformal geodesics, in the system which extends across the horizon, simplifies to
\[\lambda=\frac12V^2\dot{v}=\frac12e^\chi,\;\;\nu=\dot{u}=e^{-\chi},\;\;a_L=\frac12ae^\chi,\;\;a_N=-ae^{-\chi},\]
with constant $a$, $A=a^2$ and 
\[\dot{\chi}=a-VV'e^{-\chi},\]
which integrates to give
\[e^\chi+V^2e^{-\chi}=2ar+C\]
with constant of integration $C$. If $a\neq 0$ it is convenient to set $C=-2ar_0$ when we can derive the following equation for $r$:
\[\dot{r}^2=-V^2+a^2(r-r_0)^2=a^2(r-r_0)^2-1+\frac{2m}{r},\]
from which we obtain information about turning points. Note that 
$\dot{r}$ is zero where $V^2=a^2(r-r_0)^2$ which can be solved graphically. 
If $r_0>2m$ then there are two solutions to this equation, both greater than $2m$; for one the conformal geodesic comes in from past null infinity $\scri^-$, has a turning point and goes out to future null infinity 
$\scri^+$; for the other the conformal geodesic comes from the past 
singularity in the Penrose diagram, crosses $r=2m$ and reaches a maximum value of $r$, before recrossing $r=2m$ and ending at the future singularity. Reducing $r_0$, these solutions will persist until a critical value 
of $r_0$ is reached, when the graphs of $V^2$ and $a^2(r-r_0)^2$ are tangent; now there is a trajectory of the Killing vector $\partial/\partial t$ which is a conformal geodesic (as is to be expected from Theorem 3.1), 
and there are incoming and outgoing radial conformal 
geodesics asymptotic to it. Reducing $r_0$ still further, $\dot{r}$ has no zero and we obtain incoming or outgoing radial conformal geodesics with no turning point. 

\medskip

Another integrable case is that of solutions confined to the equatorial plane. Set $\alpha=\gamma=a_0=a_2=\cos\theta=0$. There are circles like this with $\beta=0$ too which are trajectories of the $\partial/\partial\phi$ 
Killing vector, but more generally we can set
\[\beta=\cos\chi,\;\;\;\delta=\sin\chi,\;\;a_1=-q\sin\chi,\;\;\;a_3=q\cos\chi,\]
then $q=\dot\chi+(V/r)\sin\chi=$ constant. This implies
\[(r\sin\chi)^.=qr\cos\chi=qr\dot{r}/V,\]
which is solvable as
\[\dot{r}^2=V^2(1-(G+C)^2/r^2)\]
with constant $C$ and $G'=qr/V$. The integral could be done explicitly, and there will be solutions which do not wind around the centre as well as ones which do, but the curves will not close or spiral.

%

\end{document}